\documentclass[11pt]{article}
\usepackage{amsfonts}
\usepackage{mathrsfs}
\usepackage{latexsym,amsmath}
\usepackage{amssymb,array}
\usepackage{enumitem}
\setlist{noitemsep}
\makeatletter \oddsidemargin 0in \evensidemargin 0in \textwidth 16cm
\RequirePackage[dvips]{graphicx} \textheight 20cm
\begin{document}
    \newcommand{\ov}{\overline}
    \newcommand{\om}{\omega}
    \newcommand{\ga}{\gamma}
    \newcommand{\cd}{\circledast}
    \newtheorem{thm}{Theorem}[section]
    \newtheorem{remark}{Remark}[section]
    \newtheorem{counterexample}{Counterexample}[section]
    \newtheorem{coro}{Corollary}[section]
    \newtheorem{propo}{Proposition}[section]
    \newtheorem{definition}{Definition}[section]
    \newtheorem{example}{Example}[section]
    \newtheorem{lem}{Lemma}[section]
    \numberwithin{equation}{section}
    \date{}
\title{Reliability study of proportional odds family of discrete distributions}
\author{Pradip Kundu and Asok K. Nanda\footnote
    {Corresponding author, e-mail:
    asok.k.nanda@gmail.com; asok@iiserkol.ac.in}\\
    Department of Mathematics and Statistics, IISER
Kolkata\\ Mohanpur 741246, India} \maketitle
\begin{abstract}
The proportional odds model gives a method of generating new family
of distributions by adding a parameter, called tilt parameter, to
expand an existing family of distributions. The new family of
distributions so obtained is known as Marshall-Olkin family of
distributions or Marshall-Olkin extended distributions. In this
paper, we consider Marshall-Olkin family of distributions in
discrete case with fixed tilt parameter. We study different ageing
properties, as well as different stochastic orderings of this family
of distributions. All the results of this paper are supported by
several examples.
\end{abstract}
\section{Introduction}
Since the work of Bennett \cite{benn} and Pettitt \cite{pett} on the
proportional odds (PO) model in the survival analysis context, it
has been used by many researchers. The assumption of constant hazard
ratio is unreasonable in many practical cases as discussed by
Bennett \cite{benn}, Kirmani and Gupta \cite{kirmani} and Rossini
and Tsiatis \cite{ross}. The PO model has been used by Bennett
\cite{benn} to demonstrate the effectiveness of a cure, when the
mortality rate of a group having some disease approaches that of a
(disease-free) control group as time progresses. After Bennet's
\cite{benn} work, the PO model has found many practical
applications, see, for instance, Collett \cite{coll}, Dinse and
Lagakos \cite{dinse}, Pettitt \cite{pett} and Rossini and Tsiatis
\cite{ross}. Let $X$ and $Y$ be two random variables with
distribution functions $F(\cdot)$, $G(\cdot)$, survival functions
$\bar{F}(\cdot)$, $\bar{G}(\cdot)$, probability density functions
$f(\cdot)$, $g(\cdot)$ and hazard rate functions
$r_{X}(\cdot)=f(\cdot)/\bar{F}(\cdot)$,
$r_{Y}(\cdot)=g(\cdot)/\bar{G}(\cdot)$. Let the odds functions of
$X$ and $Y$ be denoted respectively by
$\theta_{X}(t)=\bar{F}(t)/F(t)$ and $\theta_{Y}(t)=\bar{G}(t)/G(t)$.
The random variables $X$ and $Y$ are said to satisfy PO model with
proportionality constant $\alpha$ if $\theta_{Y}(t)=\alpha
\theta_{X}(t)$. For more discussion on PO models one may refer to
Kirmani and Gupta \cite{kirmani}. It is observed that, in terms of
survival functions, the PO model can be represented as
\begin{equation}\label{poalt}\bar{G}(t)=\frac{\alpha
\bar{F}(t)}{1-\bar{\alpha}\bar{F}(t)},\end{equation} where
$\bar{\alpha}=1-\alpha$. From the above representation we have
$$\frac{r_{Y}(t)}{r_{X}(t)}=\frac{1}{1-\bar{\alpha}\bar{F}(t)}=\frac{G(t)}{F(t)},$$
so that the hazard ratio is increasing (resp. decreasing) for
$\alpha>1$ (resp. $\alpha<1$) and it convergence to 1 as $t$ tends
to $\infty$. This property of hazard functions makes the PO model
reasonable in many practical applications. This is in contrast to
the proportional hazards model where the ratio of the hazard rates
remains constant with time.\par The model (\ref{poalt}), with
$0<\alpha <\infty$; gives a method of introducing new parameter
$\alpha$ to a family of distributions for obtaining more flexible
new family of distributions as discussed by Marshall and Olkin
\cite{marsh1}. The family of distributions so obtained is known as
Marshall-Olkin family of distributions or Marshall-Olkin extended
distributions (for details see \cite{marsh1,marsh2}). For more
discussion and applications of Marshall-Olkin family of
distributions one can see \cite{caron,cord1,cord,ghit2}. The
parameter $\alpha$ is called `tilt parameter'. This is because the
hazard rate of the new family is shifted below or above the hazard
rate of the underlying (baseline) distribution for $\alpha\geq 1$
and $0<\alpha \leq 1$ respectively. Thus, Marshall-Olkin family of
distributions has implications both in terms of PO model as well as
in generating a new family of flexible distributions, and hence it
is worth to investigate this family of distributions.\par Kirmani
and Gupta \cite{kirmani} studied some ageing properties of the PO
model with fixed tilt parameter $\alpha$. Nanda and Das
\cite{nanda1} studied different ageing classes of Marshall-Olkin
family of distributions taking the tilt parameter as random
variable. Ghitany and Kotz \cite{ghit} studied the reliability
properties by taking $\bar{F}$ as the reliability function of the
linear failure-rate distribution. Gupta et al. \cite{guptar}
compared the Marshall-Olkin extended distribution and the original
distribution with respect to some stochastic orderings. Nanda and
Das \cite{nanda2} compared this family of distributions with respect
to different stochastic orderings by taking the tilt parameter
random. All the studies mentioned above consider the original
(baseline) distribution to be continuous. However, not much work is
available in the literature for discrete case. In this paper, we
study different ageing properties, as well as different stochastic
orderings of this family of distributions with discrete baseline
distribution and with fixed tilt parameter.
\section{Preliminaries}
Here we discuss the survival function, the hazard (failure) rate
function, and the mean residual life of a discrete random variable
$X$ with support $\mathbb{N}=\{1,2,...\}$. Let the probability mass
function (pmf) of $X$ be given by $f(k)=P\{X=k\}$, and the
distribution function be $F(\cdot)$ so that the reliability
(survival) function $\bar{F}(\cdot)$ of X becomes
$$\bar{F}(k)=P\{X>k\}=\sum_{j=k+1}^\infty f(j),~k=1,2,...,$$ with $\bar{F}(0)=1$. The
failure rate function $r(\cdot)$ (Shaked et al. \cite{shaked}) is
given by $$r(k)=P\{X=k|X\geq k\}=\frac{P\{X=k\}}{P\{X\geq
k\}}=\frac{f(k)}{\bar{F}(k-1)},$$ and the reversed hazard rate
function is given by $\tilde{r}(k)=f(k)/F(k)$. Below we give the
definitions of different discrete ageing classes.
\begin{definition}
A discrete random variable $X$ is said to be
\begin{enumerate}[label=(\roman*)]
\item ILR (DLR) i.e. increasing (decreasing) in likelihood ratio if $f(k)$ is log-concave (log-convex),
i.e. if $f(k+2)f(k)\leq(\geq)~f^2(k+1)$, $k\in \mathbb{N}$ (Dewan
and Sudheesh \cite{dewan});
\item IFR (DFR) i.e. increasing (decreasing) failure rate if $r(k)$ is increasing (decreasing) in
$k\in \mathbb{N}$. This is equivalent to the fact that $\bar{F}(k +
1)/\bar{F}(k)$ is decreasing (increasing) in $k\in \mathbb{N}$
(Salvia and Bollinger \cite{salvia}, Shaked et al. \cite{shaked},
Gupta et al. \cite{gupta});
\item IFRA (DFRA) i.e. increasing (decreasing) in failure rate average if $[\bar{F}(k)]^{1/k}$ is decreasing (increasing) in $k$, i.e., $[\bar{F}(k)]^{1/k}\geq(\leq)~[\bar{F}(k+1)]^{1/(k+1)}$, $k\in \mathbb{N}$
(Esary et al. \cite{esary}, Shaked et al. \cite{shaked});
\item NBU (NWU) i.e. new better (worse) than used if
$\bar{F}(j + k)\leq(\geq)~\bar{F}(j)\bar{F}(k),~j,k\in \mathbb{N}$
(Esary et al. \cite{esary}, Shaked et al. \cite{shaked});
\item DRHR (decreasing reversed hazard rate) if $\tilde{r}(k)$ is decreasing in
$k$, i.e. if $F(k)$ is log-concave, i.e. if $[F(k+1)]^2\geq F(k)
F(k+2)$, $k\in \mathbb{N}$ (Nanda and Sengupta \cite{nanda}, Li and
Xu \cite{lix}).
\item NBAFR (new better than used in failure rate average) if
$[\bar{F}(k)]^{1/k}\leq \bar{F}(1)$, $k\in \mathbb{N}$ (Fagiuoli and
Pellerey \cite{fagi}).
\end{enumerate}
\end{definition}
\subsection{Proportional odds family of discrete distributions} Let $X$
be a discrete random variable with support $\mathbb{N}=\{1,2,...\}$
having pmf $f(\cdot)$, distribution function $F(\cdot)$, survival
function $\bar{F}(\cdot)$, hazard rate function $r_X(\cdot)$, and
reversed hazard rate function $\tilde{r}_X(\cdot)$. Starting with
the survival function $\bar{F}$, the survival function of the
proportional odds family (also known as Marshall-Olkin family) of
discrete distribution is given by
\begin{equation}\label{eqmos}\bar{G}(k;\alpha)=\frac{\alpha \bar{F}(k)}{1-\bar{\alpha}\bar{F}(k)},~k=1,2,...,~0<\alpha <\infty,~\bar{\alpha}=1-\alpha,\end{equation}
with $\bar{G}(0;\alpha)=1$. Let the corresponding random variable be
denoted by $Y$. Now the distribution function of $Y$ is given by
\begin{equation}\label{eqmod}G(k;\alpha)=1-\bar{G}(k;\alpha)=\frac{F(k)}{1-\bar{\alpha}\bar{F}(k)},\end{equation}
whereas the pmf is given by
\begin{equation}\label{eqmom}g(k;\alpha)=\bar{G}(k-1;\alpha)-\bar{G}(k;\alpha)=\frac{\alpha f(k)}{[1-\bar{\alpha}\bar{F}(k-1)][1-\bar{\alpha}\bar{F}(k)]}.\end{equation}
The corresponding hazard rate and the reversed hazard rate functions
are given by
\begin{equation}\label{eqmoh}r_{Y}(k;\alpha)=\frac{r_X(k)}{1-\bar{\alpha}\bar{F}(k)},\end{equation}
\begin{equation}\label{eqmorh}\tilde{r}_{Y}(k;\alpha)=\frac{\alpha \tilde{r}_X(k)}{1-\bar{\alpha}\bar{F}(k-1)}.\end{equation}
It is to be mentioned here that different properties of
(\ref{eqmos}) have been studied by D\'{e}niz and Sarabia
\cite{deniz} by taking $F$ as the cdf of Poisson random variable.
\section{Stochastic Ageing properties}
In this section we study how different ageing properties of $X$ are
transmitted to the random variable $Y$.\par With the following two
counterexamples, one with $\alpha>1$ and the other with $\alpha<1$,
we show that if $X$ is ILR, then $Y$ is neither ILR nor DLR.
\begin{counterexample}
Consider the random variable $X$ with the mass function given by
$$f(k)=\left\{
                                          \begin{array}{ll}
                                            0, & \hbox{if $k=1$;} \\
                                            0.1, & \hbox{if $k=2$;} \\
                                            0.25, & \hbox{if $k=3$;} \\
                                            0.35, & \hbox{if $k=4$;} \\
                                            0.3, & \hbox{if $k=5$.}
                                          \end{array}\right.$$
Clearly $X$ is ILR. For $\alpha=5$, we have the mass function of $Y$
as
$$g(k;5)=\left\{
                                          \begin{array}{ll}
                                            0, & \hbox{if $k=1$;} \\
                                            \frac{1}{46}, & \hbox{if $k=2$;} \\
                                            \frac{125}{1656}, & \hbox{if $k=3$;} \\
                                            \frac{175}{792}, & \hbox{if $k=4$;} \\
                                            \frac{15}{22}, & \hbox{if $k=5$.}
                                          \end{array}
                                        \right.$$ It is observed
that $Y$ is neither ILR nor DLR.
\end{counterexample}
\begin{counterexample}
Consider the random variable $X$ with mass function given by
$$f(k)=\left\{
                                          \begin{array}{ll}
                                            0, & \hbox{if $k=1$;} \\
                                            0.3, & \hbox{if $k=2$;} \\
                                            0.34, & \hbox{if $k=3$;} \\
                                            0.26, & \hbox{if $k=4$;} \\
                                            0.1, & \hbox{if $k=5$.}
                                          \end{array}\right.$$ Here
$X$ is ILR. For $\alpha=0.2$, we have the mass function of $Y$ as
$$g(k;0.2)=\left\{
                                          \begin{array}{ll}
                                            0, & \hbox{if $k=1$;} \\
                                            \frac{15}{22}, & \hbox{if $k=2$;} \\
                                            \frac{425}{1958}, & \hbox{if $k=3$;} \\
                                            \frac{325}{4094}, & \hbox{if $k=4$;} \\
                                            \frac{1}{46}, & \hbox{if $k=5$.}
                                          \end{array}
                                        \right.$$ It is observed
that $Y$ is neither ILR nor DLR.$\hfill\Box$
\end{counterexample}
\par With the following two counterexamples, one for $\alpha>1$ and
the other for $\alpha<1$, we show that if $X$ is DLR, then $Y$ is
neither DLR nor ILR.
\begin{counterexample}
Consider the random variable $X$ with mass function given by
$$f(k)=\left\{
                                          \begin{array}{ll}
                                            0.36, & \hbox{if $k=1$;} \\
                                            0.26, & \hbox{if $k=2$;} \\
                                            0.21, & \hbox{if $k=3$;} \\
                                            0.17, & \hbox{if $k=4$.}
                                          \end{array}\right.$$ Here
$X$ is DLR. For $\alpha=2$, we have the mass function of $Y$ as
$$g(k;2)=\left\{
                                          \begin{array}{ll}
                                            \frac{9}{41}, & \hbox{if $k=1$;} \\
                                            \frac{650}{2829}, & \hbox{if $k=2$;} \\
                                            \frac{700}{2691}, & \hbox{if $k=3$;} \\
                                            \frac{34}{117}, & \hbox{if $k=4$.}
                                          \end{array}
                                        \right.$$ It is observed
that $Y$ is neither DLR nor ILR.
\end{counterexample}
\begin{counterexample}
Consider the random variable $X$ with mass function given by
$$f(k)=\left\{
                                          \begin{array}{ll}
                                            0.26, & \hbox{if $k=1$;} \\
                                            0.18, & \hbox{if $k=2$;} \\
                                            0.24, & \hbox{if $k=3$;} \\
                                            0.32, & \hbox{if $k=4$.}
                                          \end{array}\right.$$ Here
$X$ is DLR. For $\alpha=0.4$, we have the mass function of $Y$ as
$$g(k;0.4)=\left\{
                                          \begin{array}{ll}
                                            \frac{65}{139}, & \hbox{if $k=1$;} \\
                                            \frac{2250}{11537}, & \hbox{if $k=2$;} \\
                                            \frac{1500}{8383}, & \hbox{if $k=3$;} \\
                                            \frac{16}{101}, & \hbox{if $k=4$.}
                                          \end{array}
                                        \right.$$ It is observed
that $Y$ is neither DLR nor ILR.$\hfill\Box$
\end{counterexample}
\par The following theorem gives the condition under which IFR/DFR
property of $X$ is transmitted to the random variable $Y$. The proof
follows from the fact that $1/\left(1-\bar{\alpha}\bar{F}(k)\right)$
is increasing (resp. decreasing) in $k$ for $\alpha\geq (resp. \leq)
1$.
\begin{thm}If $X$ is IFR (resp. DFR) and $\alpha\geq (resp. \leq)~1$, then $Y$ is IFR (resp. DFR).$\hfill\Box$\end{thm}
\par The following counterexample shows that if $\alpha<1$, then the
IFR property of $X$ may not be transmitted to the random variable
$Y$.
\begin{counterexample}\label{exdifr}
Consider the random variable $X$ following discrete IFR distribution
(cf. Salvia and Bollinger \cite{salvia}) with
$$f(k)=(k-c)c^{k-1}/k!,~k\in \mathbb{N},~0<c\leq 1.$$ Here
$\bar{F}(k)=c^{k}/k!$ and $r_X(k)=1-c/k$ so that $X$ is IFR. Now
$$r_{Y}(k;\alpha)=\frac{1-c/k}{1-\bar{\alpha}c^{k}/k!}.$$
It is observed that, for $c=0.8$ and $\alpha=0.2$, we have
$r_{Y}(2;0.2)=0.8064516$, $r_{Y}(3;0.2)=0.7870635$, and
$r_{Y}(4;0.2)=0.8110739$. This shows that $Y$ is neither IFR nor
DFR. $\hfill\Box$\end{counterexample} \par The following
counterexample shows that the DFR property of $X$ may not be
transmitted to the random variable $Y$ when $\alpha>1$.
\begin{counterexample}\label{ext1wd}
Let $X$ follow the Type I discrete Weibull distribution (cf.
Nakagawa and Osaki \cite{nakag}) with pmf given by
$$f(k)=q^{(k-1)^\beta}-q^{k^\beta},~k\in \mathbb{N},~q\in(0,1),~\beta>0.$$ Then the corresponding survival function is given by
$\bar{F}(k)=q^{k^\beta}$, and the hazard rate function is given by
$r_X(k)=1-q^{k^\beta -(k-1)^\beta}$. Here $X$ is DFR for
$0<\beta<1$. Note that
$$r_{Y}(k;\alpha)=\frac{1-q^{k^\beta -(k-1)^\beta}}{1-\bar{\alpha}q^{k^\beta}}.$$
It is observed that, for $\beta=0.8$, $\alpha=5$ and $q=0.5$, we
have $r_{Y}(7;5)=0.2759209$, $r_{Y}(10;5)=0.2834942$, and
$r_{Y}(13;5)=0.2793229$. This shows that $Y$ is neither DFR nor
IFR.$\hfill\Box$
\end{counterexample}
\par Below we see that, for $\alpha\geq (resp. \leq)~1$, the NBU
(resp. NWU) property of $X$ is transmitted to the random variable
$Y$.
\begin{thm}
If $X$ is NBU (resp. NWU) and $\alpha\geq (resp. \leq)~1$, then $Y$
is NBU (resp. NWU).\end{thm} \textbf{Proof:} Let $X$ be NBU (resp.
NWU). Then $Y$ will be  NBU (resp. NWU) if and only if
$$\frac{1-\bar{\alpha} \bar{F}(j+k)}{\bar{F}(j+k)}\geq (resp. \leq)
\frac{(1-\bar{\alpha}\bar{F}(j))(1-\bar{\alpha}\bar{F}(k))}{\alpha
\bar{F}(j)\bar{F}(k)}.$$ This is equivalent to the fact that
$$\frac{1}{\bar{F}(j+k)} \geq (resp. \leq)
\frac{1-\bar{\alpha}\left(\bar{F}(j)+\bar{F}(k)\right)+\bar{\alpha}\bar{F}(k)\bar{F}(j)}{\alpha
\bar{F}(j)\bar{F}(k)},$$ which holds if $$\frac{1}{\bar{F}(j+k)}
\geq (resp. \leq) \frac{1}{\bar{F}(j)\bar{F}(k)}.$$ The last
inequality follows from the fact that, for $\alpha\geq (resp.
\leq)~1$,
$$1-\bar{\alpha}\left(\bar{F}(j)+\bar{F}(k)\right)+\bar{\alpha}\bar{F}(k)\bar{F}(j)\leq (resp. \geq)~
\alpha.$$ Hence the theorem follows.$\hfill\Box$
\par The following counterexamples show that, for $\alpha< (resp.
>)~1$, the NBU (resp. NWU) property of $X$ may not be transmitted to
the random variable $Y$.
\begin{counterexample}\label{notnbu}
Consider $X$ following the discrete S-distribution (cf. Bracuemond
and Gaudoin \cite{brac}) with pmf given by
$$f(k)=p(1-a^k)\prod_{i=1}^{k-1}(1-p+p a^i),~k\in \mathbb{N},~
0<p\leq 1,~0<a<1.$$ This gives the survival function as
$\bar{F}(k)=\prod_{i=1}^{k}(1-p+p a^i)$, and hazard rate function as
$r_X(k)=p(1-a^k)$. Here $X$ is NBU. Now
$$\bar{G}(k;\alpha)=\frac{\alpha \prod_{i=1}^{k}(1-p+p a^i)}{1-\bar{\alpha}\prod_{i=1}^{k}(1-p+p a^i)}.$$
For $j=2$, $k=3$, $p=0.3$, $a=0.6$, $\alpha=0.2$, we have
$\bar{G}(j+k;\alpha)=0.075737$ and
$\bar{G}(j;\alpha)\bar{G}(k;\alpha)=0.063494$. This shows that $Y$
is not NBU.
\end{counterexample}
\begin{counterexample}
Let $X$ follow the distribution as given in Counterexample
\ref{ext1wd}. Then clearly $X$ is NWU for $\beta\in(0,1)$. Now, for
$j=2$, $k=3$, $\alpha=5$, $q=0.5$, we have
$\bar{G}(j+k;\alpha)=0.3062174$ and
$\bar{G}(j;\alpha)\bar{G}(k;\alpha)=0.3657684$. This shows that $Y$
is not NWU.$\hfill\Box$
\end{counterexample}
\par Kirmani and Gupta \cite{kirmani} have observed that if $X$ is
IFRA (DFRA), then $Y$ is IFRA (DFRA) for $\alpha>(<)~1$. Below we
show that if $X$ is IFRA, then $Y$ may not be IFRA or DFRA when
$\alpha<1$.
\begin{counterexample}
Let $X$ follow the distribution as given in Counterexample
\ref{notnbu}. Here $X$ is IFRA. Now, for $\alpha=0.2$, $p=0.5$, and
$a=0.6$, we have
$$[\bar{G}(k;0.2)]^{1/k}=\left\{
  \begin{array}{ll}
    0.44444, & \hbox{for $k=1$;} \\
    0.438901, & \hbox{for $k=2$;} \\
    0.457806, & \hbox{for $k=4$.}
  \end{array}
\right.$$ This shows that $Y$ is neither IFRA nor DFRA.$\hfill\Box$
\end{counterexample}
\par Below we show that if $X$ is DFRA, then, for $\alpha>1$, $Y$ may
not be DFRA or IFRA.
\begin{counterexample}
Let $X$ follow the discrete Pareto distribution with survival
function $$\bar{F}(k)=\left(\frac{d}{k+d}\right)^c,~k\in
\mathbb{N},~c,d>0,$$ which is DFR and hence DFRA. Now
$$\bar{G}(k;\alpha)=\frac{\alpha \left(\frac{d}{k+d}\right)^c}{1-\bar{\alpha}\left(\frac{d}{k+d}\right)^c}.$$
For $\alpha=6$, $d=2$, $c=3$, we have
$$[\bar{G}(k;6)]^{1/k}=\left\{
    \begin{array}{ll}
      0.7164179, & \hbox{for $k=1$;} \\
      0.658037, & \hbox{for $k=4$;} \\
      0.68081, & \hbox{for $k=8$,}
    \end{array}
  \right.$$
which is neither increasing nor decreasing in $k$, i.e. $Y$ is
neither DFRA nor IFRA.$\hfill\Box$
\end{counterexample}
\par Following theorem shows that, for $\alpha\leq 1$, the DRHR
property of $X$ is transmitted to the random variable $Y$. The proof
follows from the fact that
$1/\left(1-\bar{\alpha}\bar{F}(k-1)\right)$ is decreasing in $k$,
for $\alpha\leq 1$.
\begin{thm}If $X$ is DRHR, then $Y$ is DRHR for $\alpha\leq 1$.$\hfill\Box$\end{thm}
\par The following counterexample shows that, for $\alpha>1$, DRHR
property of $X$ may not be transmitted to the random variable $Y$.
\begin{counterexample}
Consider the random variable $X$ having distribution function given
by
$$F(k)=\left\{
                                          \begin{array}{ll}
                                            0, & \hbox{if $1\leq k<2$;} \\
                                            \frac{4}{25}, & \hbox{if $2\leq k<3$;} \\
                                            \frac{2}{5}, & \hbox{if $3\leq k<4$;} \\
                                            \frac{2}{3}, & \hbox{if $4\leq k<5$;} \\
                                            1, & \hbox{if $k\geq 5$.}
                                          \end{array}
                                        \right.$$ Clearly $X$ is DRHR. For $\alpha=4$, the distribution function of $Y$ is given by
$$G(k;4)=\left\{
                                          \begin{array}{ll}
                                            0, & \hbox{if $1\leq k<2$;} \\
                                            \frac{1}{22}, & \hbox{if $2\leq k<3$;} \\
                                            \frac{1}{7}, & \hbox{if $3\leq k<4$;} \\
                                            \frac{1}{3}, & \hbox{if $4\leq k<5$;} \\
                                            1, & \hbox{if $k\geq 5$,}
                                          \end{array}
                                        \right.$$ which is not DRHR.$\hfill\Box$
\end{counterexample}
\par The following counterexample shows that, for $\alpha<1$, NBAFR
property of $X$ may not be transmitted to the random variable $Y$.
\begin{counterexample}
Consider the random variable $X$ with reliability function given by
$$\bar{F}(k)=\left\{
                                          \begin{array}{ll}
                                             \frac{4}{5}, & \hbox{if $1\leq k<2$;} \\
                                            \frac{8}{13}, & \hbox{if $2\leq k<3$;} \\
                                            \frac{1}{2}, & \hbox{if $3\leq k<4$;} \\
                                            0, & \hbox{if $k\geq 4$.}
                                          \end{array}
                                        \right.$$
Here $X$ is NBAFR. For $\alpha=0.4$, we have the reliability
function of $Y$ as
$$\bar{G}(k;0.4)=\left\{
                                          \begin{array}{ll}
                                            \frac{8}{13}, & \hbox{if $1\leq k<2$;} \\
                                            \frac{16}{41}, & \hbox{if $2\leq k<3$;} \\
                                            \frac{2}{7}, & \hbox{if $3\leq k<4$;} \\
                                            0, & \hbox{if $k\geq 4$.}
                                          \end{array}
                                        \right.$$
It is clear that $Y$ is not NBAFR.$\hfill\Box$\end{counterexample}
We summarize the above findings in Table 1.
\begin{table}[h]
\begin{center} Table 1: Preservation of ageing classes\\
\begin{tabular}{|l|l|l|}
\hline
Ageing properties  & \multicolumn{1}{c|}{$\alpha<1$} &
\multicolumn{1}{c|}{$\alpha>1$}
\\\hline ILR & Not Preserved & Not Preserved
 \\
DLR & Not Preserved & Not Preserved
 \\
IFR & Not Preserved & Preserved
 \\
DFR & Preserved & Not Preserved
 \\
NBU & Not Preserved & Preserved
 \\
NWU & Preserved & Not Preserved
 \\
IFRA & Not Preserved & Preserved
 \\
DFRA & Preserved & Not Preserved
 \\
DRHR & Preserved & Not Preserved
 \\
NBAFR & Not Preserved &
 \\
\hline
\end{tabular}
\end{center}
\end{table}
\section{Stochastic Orderings}
Let $X_{1}$ and $X_{2}$ be two discrete random variables with
support $\mathbb{N}=\{1,2,...\}$ having respective pmf
$f_{1}(\cdot)$, $f_{2}(\cdot)$, distribution function
$F_{1}(\cdot)$, $F_{2}(\cdot)$, and survival function
$\bar{F}_{1}(\cdot)$, $\bar{F}_{2}(\cdot)$. Let the survival
function of the Marshall-Olkin family of discrete distributions be
given by $$\bar{G}_{i}(k;\alpha)=\frac{\alpha
\bar{F}_{i}(k)}{1-\bar{\alpha} \bar{F}_{i}(k)},~0<\alpha<\infty,~
\bar{\alpha}=1-\alpha,$$ and let the corresponding random variable
be $Y_{i}$, $i=1,2$. The following theorem shows that the usual
stochastic order between $X_{1}$ and $X_{2}$ and that of $Y_{1}$ and
$Y_{2}$ are equivalent.
\begin{thm}\label{thst}
$Y_1\leq_{st} Y_2$ if and only if $X_1\leq_{st} X_2$.
\end{thm}
\textbf{Proof:} Note that $Y_1\leq_{st} Y_2$ if, and only if
\begin{eqnarray*}
\frac{\alpha\bar{F}_{1}(k)}{1-\bar{\alpha}\bar{F}_{1}(k)}&\leq&
\frac{\alpha\bar{F}_{2}(k)}{1-\bar{\alpha}\bar{F}_{2}(k)},
\end{eqnarray*} which is equivalent to the fact that $\bar{F}_{1}(k)\leq
\bar{F}_{2}(k)$. Hence the theorem follows.$\hfill\Box$ \par The
following theorem gives condition on $\alpha$, under which hazard
rate order between $X_1$ and $X_2$ is transmitted to that between
$Y_1$ and $Y_1$.
\begin{thm}\label{thhr}
If $X_1\leq_{hr} X_2$, then $Y_1\leq_{hr} Y_2$, provided $\alpha\geq
1$.
\end{thm}
\textbf{Proof:} Since hazard rate order is stronger than usual
stochastic order, we have, for $\alpha\geq 1$,
$$\frac{1}{1-\bar{\alpha}\bar{F}_{1}(k)}\geq\frac{1}{1-\bar{\alpha}\bar{F}_{2}(k)}.$$
Now, using the hypothesis we have, from (\ref{eqmoh}),
$$r_{Y_{1}}(k;\alpha)=\frac{r_{X_{1}}(k)}{1-\bar{\alpha}
\bar{F}_{1}(k)}\geq \frac{r_{X_{2}}(k)}{1-\bar{\alpha}
\bar{F}_{2}(k)}=r_{Y_{2}}(k;\alpha).$$ Hence the theorem follows.$\hfill\Box$\\
\par The following counterexample shows that the above theorem does
not hold if $\alpha<1$.
\begin{counterexample}
Consider the random variables $X_1$ and $X_2$ with respective
reliability function $$\bar{F}_{1}(k)=\left\{
                                          \begin{array}{ll}
                                            1, & \hbox{if $1\leq k<2$;} \\
                                            \frac{1}{2}, & \hbox{if $2\leq k<3$;} \\
                                            \frac{2}{5}, & \hbox{if $3\leq k<4$;} \\
                                            0, & \hbox{if $k\geq 4$,}
                                          \end{array}
                                        \right.$$
and
$$\bar{F}_{2}(k)=\left\{
                                          \begin{array}{ll}
                                            1, & \hbox{if $1\leq k<2$;} \\
                                            \frac{5}{8}, & \hbox{if $2\leq k<3$;} \\
                                            \frac{11}{20}, & \hbox{if $3\leq k<4$;} \\
                                            0, & \hbox{if $k\geq 4$.}
                                          \end{array}
                                        \right.$$
This shows that $X_1\leq_{hr} X_2$. For $\alpha=0.2$, we have the
respective reliability function of $Y_1$ and $Y_2$ as
$$\bar{G}_{1}(k;\alpha)=\left\{
                                          \begin{array}{ll}
                                            1, & \hbox{if $1\leq k<2$;} \\
                                            \frac{1}{6}, & \hbox{if $2\leq k<3$;} \\
                                            \frac{2}{17}, & \hbox{if $3\leq k<4$;} \\
                                            0, & \hbox{if $k\geq 4$,}
                                          \end{array}
                                        \right.$$
and
$$\bar{G}_{2}(k;\alpha)=\left\{
                                          \begin{array}{ll}
                                            1, & \hbox{if $1\leq k<2$;} \\
                                            \frac{1}{4}, & \hbox{if $2\leq k<3$;} \\
                                            \frac{11}{56}, & \hbox{if $3\leq k<4$;} \\
                                            0, & \hbox{if $k\geq 4$.}
                                          \end{array}
                                        \right.$$ This shows that $Y_1\nleq_{hr} Y_2$.$\hfill\Box$
\end{counterexample}
\par The following theorem gives the condition on $\alpha$ such that
reversed hazard rate order between  $X_1$ and $X_2$ is transmitted
to that between $Y_1$ and $Y_2$.
\begin{thm}\label{thrhr}
If $X_1\leq_{rhr} X_2$, then $Y_1\leq_{rhr} Y_2$, for $0<\alpha\leq
1$.
\end{thm}
\textbf{Proof:} Since reversed hazard rate order is stronger than
usual stochastic order, we have, for $\alpha\leq 1$,
$$\frac{1}{1-\bar{\alpha}\bar{F}_{1}(k-1)}\leq\frac{1}{1-\bar{\alpha}\bar{F}_{2}(k-1)}.$$
Now, using the hypothesis we have, from (\ref{eqmorh}),
$$\tilde{r}_{Y_1}(k;\alpha)=\frac{\alpha
\tilde{r}_{X_{1}}(k)}{1-\bar{\alpha}\bar{F_1}(k-1)}\leq \frac{\alpha
\tilde{r}_{X_{2}}(k)}{1-\bar{\alpha}\bar{F_2}(k-1)}=\tilde{r}_{Y_2}(k;\alpha).$$
Hence the theorem follows.$\hfill\Box$\par That the above theorem
does not hold in case of $\alpha>1$ is shown in the following
counterexample.
\begin{counterexample}
Consider the random variables $X_1$ and $X_2$ with respective
distribution function $$F_{1}(k)=\left\{
                                          \begin{array}{ll}
                                            0, & \hbox{if $1\leq k<2$;} \\
                                            \frac{5}{24}, & \hbox{if $2\leq k<3$;} \\
                                            \frac{1}{2}, & \hbox{if $3\leq k<4$;} \\
                                            \frac{3}{4}, & \hbox{if $4\leq k<5$;} \\
                                            1, & \hbox{if $k\geq 5$,}
                                          \end{array}
                                        \right.$$ and $$F_{2}(k)=\left\{
                                          \begin{array}{ll}
                                           0, & \hbox{if $1\leq k<2$;} \\
                                            \frac{1}{6}, & \hbox{if $2\leq k<3$;} \\
                                            \frac{5}{12}, & \hbox{if $3\leq k<4$;} \\
                                            \frac{2}{3}, & \hbox{if $4\leq k<5$;} \\
                                            1, & \hbox{if $k\geq 5$.}
                                          \end{array}
                                        \right.$$
Clearly $X_1\leq_{rhr} X_2$. For $\alpha=4$, we have the
distribution functions of $Y_1$ and $Y_2$ respectively as
$$G_{1}(k;4)=\left\{
                                          \begin{array}{ll}
                                            0, & \hbox{if $1\leq k<2$;} \\
                                            \frac{5}{81}, & \hbox{if $2\leq k<3$;} \\
                                            \frac{1}{5}, & \hbox{if $3\leq k<4$;} \\
                                            \frac{3}{7}, & \hbox{if $4\leq k<5$;} \\
                                            1, & \hbox{if $k\geq 5$,}
                                          \end{array}
                                        \right.$$
and
$$G_{2}(k;4)=\left\{
                                          \begin{array}{ll}
                                           0, & \hbox{if $1\leq k<2$;} \\
                                            \frac{1}{21}, & \hbox{if $2\leq k<3$;} \\
                                            \frac{5}{33}, & \hbox{if $3\leq k<4$;} \\
                                            \frac{1}{3}, & \hbox{if $4\leq k<5$;} \\
                                            1, & \hbox{if $k\geq 5$.}
                                          \end{array}
                                        \right.$$ This shows
that $Y_1\nleq_{rhr} Y_2$.$\hfill\Box$
\end{counterexample}
\par Following two counterexamples show that the likelihood ratio
order between $X_1$ and $X_2$ is not necessarily transmitted to that
between $Y_1$ and $Y_2$.
\begin{counterexample}
Let $X_1$ and $X_2$ have the respective probability mass function
$$f_{1}(k)=\left\{
                                          \begin{array}{ll}
                                            0, & \hbox{if $k=1$;} \\
                                            0.3, & \hbox{if $k=2$;} \\
                                            0.4, & \hbox{if $k=3$;} \\
                                            0.2, & \hbox{if $k=4$;} \\
                                            0.1, & \hbox{if $k=5$,}
                                          \end{array}
                                        \right.$$
and
$$f_{2}(k)=\left\{
                                          \begin{array}{ll}
                                            0, & \hbox{if $k=1$;} \\
                                            0.2, & \hbox{if $k=2$;} \\
                                            0.3, & \hbox{if $k=3$;} \\
                                            0.2, & \hbox{if $k=4$;} \\
                                            0.3, & \hbox{if $k=5$.}
                                          \end{array}
                                        \right.$$
Clearly $X_1\leq_{lr} X_2$. For $\alpha=5$, we have the mass
functions of $Y_1$ and $Y_2$ respectively as
$$g_{1}(k;5)=\left\{
                                          \begin{array}{ll}
                                            0, & \hbox{if $k=1$;} \\
                                            \frac{3}{38}, & \hbox{if $k=2$;} \\
                                            \frac{50}{209}, & \hbox{if $k=3$;} \\
                                            \frac{25}{77}, & \hbox{if $k=4$;} \\
                                            \frac{5}{14}, & \hbox{if $k=5$,}
                                          \end{array}
                                        \right.$$
and
$$g_{2}(k;5)=\left\{
                                          \begin{array}{ll}
                                            0, & \hbox{if $k=1$;} \\
                                            \frac{1}{21}, & \hbox{if $k=2$;} \\
                                            \frac{5}{42}, & \hbox{if $k=3$;} \\
                                            \frac{5}{33}, & \hbox{if $k=4$;} \\
                                            \frac{15}{22}, & \hbox{if $k=5$.}
                                          \end{array}
                                        \right.$$ This shows
that $Y_1\nleq_{lr} Y_2$.$\hfill\Box$
\end{counterexample}
\begin{counterexample}
Take the random variables $X_1$ and $X_2$ having respective mass
functions
$$f_{1}(k)=\left\{
                                          \begin{array}{ll}
                                            0, & \hbox{if $k=1$;} \\
                                            0.3, & \hbox{if $k=2$;} \\
                                            0.3, & \hbox{if $k=3$;} \\
                                            0.2, & \hbox{if $k=4$;} \\
                                            0.2, & \hbox{if $k=5$,}
                                          \end{array}
                                        \right.$$
and
$$f_{2}(k)=\left\{
                                          \begin{array}{ll}
                                            0, & \hbox{if $k=1$;} \\
                                            0.2, & \hbox{if $k=2$;} \\
                                            0.3, & \hbox{if $k=3$;} \\
                                            0.24, & \hbox{if $k=4$;} \\
                                            0.26, & \hbox{if $k=5$.}
                                          \end{array}
                                        \right.$$
Clearly $X_1\leq_{lr} X_2$. For $\alpha=0.2$, we have the mass
functions of $Y_1$ and $Y_2$ as
$$g_{1}(k;0.2)=\left\{
                                          \begin{array}{ll}
                                            0, & \hbox{if $k=1$;} \\
                                            \frac{15}{22}, & \hbox{if $k=2$;} \\
                                            \frac{75}{374}, & \hbox{if $k=3$;} \\
                                            \frac{25}{357}, & \hbox{if $k=4$;} \\
                                            \frac{1}{21}, & \hbox{if $k=5$,}
                                          \end{array}
                                        \right.$$
and
$$g_{2}(k;0.2)=\left\{
                                          \begin{array}{ll}
                                            0, & \hbox{if $k=1$;} \\
                                            \frac{5}{9}, & \hbox{if $k=2$;} \\
                                            \frac{5}{18}, & \hbox{if $k=3$;} \\
                                            \frac{10}{99}, & \hbox{if $k=4$;} \\
                                            \frac{13}{198}, & \hbox{if $k=5$.}
                                          \end{array}
                                        \right.$$ This shows
that $Y_1\nleq_{lr} Y_2$.$\hfill\Box$
\end{counterexample}
We summarize the above findings in Table 2.
\begin{table}[h]
\begin{center} Table 2: Preservation of stochastic orderings\\
\begin{tabular}{|l|l|l|}
\hline
Stochastic orders between & \multicolumn{1}{c|}{$\alpha<1$} &
\multicolumn{1}{c|}{$\alpha>1$}\\baseline distributions   &
 & 
\\\hline Usual stochastic order & Preserved & Preserved
 \\
Hazard rate order & Not Preserved & Preserved
 \\
Reversed hazard rate order & Preserved & Not Preserved
 \\
Likelihood ratio order & Not Preserved & Not Preserved
 \\
\hline
\end{tabular}
\end{center}
\end{table}
\section{Conclusion}Marshall and Olkin \cite{marsh1} introduced a method of
adding a new parameter, called tilt parameter, to a family of
distributions for obtaining more flexible new families of
distributions. In the literature, some reliability properties of
this family of distributions are studied with continuous baseline
distributions. However, not much study is done in the literature, to
the best of our knowledge, for discrete baseline distributions. This
paper discusses various stochastic ageing properties, as well as
different stochastic orderings of this family with discrete baseline
distributions.
\subsection*{Acknowledgements:} The support received from IISER
Kolkata to carry out this research work is gratefully acknowledged
by Pradip Kundu. The financial support from NBHM, Govt. of India
(vide Ref. No. 2/48(25)/2014/NBHM(R.P.)/R\&D II/1393 dt. Feb. 3,
2015) is duly acknowledged by Asok K. Nanda.

\end{document}